\documentclass[12pt]{amsart}\usepackage{}\usepackage{cases}\usepackage{amsmath}\usepackage{amssymb}\usepackage{amsfonts}
\newcommand{\Pf}{\noindent {\bf Proof}}

\newtheorem{thm}{Theorem}[section]
\newtheorem{prop}{Proposition}[section]
\newtheorem{lemma}{Lemma}[section]

\newtheorem{remark}{Remark}[section]

\begin{document}
\title[GSRS half harmonic Weyl curvature] {Gradient shrinking Ricci solitons of half harmonic Weyl curvature}
\author{Jia-Yong Wu}
\address{Department of Mathematics, Shanghai Maritime University, Haigang Avenue 1550, Shanghai 201306, P. R. China.}\email{jywu81@yahoo.com}
\author{Peng Wu}
\address{Department of Mathematics, Cornell University, Ithaca, NY 14853, U.S.A.} \email{wupenguin@math.cornell.edu}
\author{William Wylie}
\address{Department of Mathematics, Syracuse University, Syracuse, NY 13244, U.S.A.} \email{wwylie@syr.edu}
\thanks{}
\dedicatory{}
\date{September 5, 2014}

\begin{abstract}
We prove that a four-dimensional gradient shrinking Ricci soliton with $\delta W^{\pm}=0$ is either Einstein, or a finite quotient of $S^3\times\mathbb{R}$, $S^2\times\mathbb{R}^2$ or $\mathbb{R}^4$. We also prove that a four-dimensional cscK gradient Ricci soliton is either K\"ahler-Einstein, or a finite quotient of $M\times\mathbb{C}$, where $M$ is a Riemann surface. The main arguments are curvature decompositions, the Weitzenb\"ock formula for half Weyl curvature, and the maximum principle.
\end{abstract}
\maketitle

\section{Introduction}
In this paper we investigate four-dimensional gradient shrinking Ricci solitons with half harmonic Weyl curvature. A Riemannian metric $g$ on a smooth manifold $M^n$ is called a gradient Ricci soliton, if there exists an $f\in C^{\infty}(M)$ and a $\lambda\in\mathbb{R}$, such that
\begin{equation} \label{solitoneq}
\begin{split}
\mathrm{Ric}+\nabla^2f=\lambda g.
\end{split}
\end{equation}
The function $f$ is called a potential function for the gradient Ricci soliton. By convention we denote it by  the triple $(M^n,g,f)$.
The gradient Ricci soliton is called shrinking, steady, or expanding, if $\lambda>0$, $\lambda=0$, or $\lambda<0$, respectively.
Gradient Ricci solitons play an important role in the Ricci flow as they are self-similar solutions to the flow,  making them possible singularity models  and critical points of Perelman's entropies.   See \cite{Cao} for an excellent survey of Ricci solitons.

Any Einstein metric is also a gradient Ricci soliton with $f$ taken to be a constant function.   The only space which is both Einstein and a gradient Ricci soliton with a non-constant potential function is a \emph{Gaussian} which is  flat Euclidean space.  If $\lambda \neq 0$, $f$ can be chosen as  $f(x) = \frac{\lambda}{2} d(p,x)^2$ and $d(p,x)$ is the Euclidean distance to some point $p \in \mathbb{R}^n$.   A simplest way to construct a non-Einstein gradient Ricci soliton is to take the product of an Einstein metric and a Gaussian.

In dimensions $2$ and $3$ these are the only examples of shrinking gradient Ricci solitons as follows from the work of  Hamilton \cite{Hamilton2}, Ivey \cite{Ivey}, Perelman \cite{Perel2}, L. Ni, Wallach \cite{NW}, and H.-D. Cao, B.-L. Chen, and X.-P. Zhu \cite{CCZ}.  In higher dimensions many authors have studied the classification of gradient Ricci solitons.  For the purposes of this introduction we focus only on results for shrinking gradient Ricci solitons.  Any locally conformal flat gradient shrinking Ricci soliton is a finite quotient of $S^n$, $S^{n-1}\times\mathbb{R}$, or $\mathbb{R}^n$ as follows from the works of various authors, see \cite{ENM, CWZ, NW, Zhang, PW, MS}. Recall that a Riemannian manifold is locally conformally flat if the Weyl tensor vanishes.  Fern\'andez-L\'opez, Garc\'ia-R\'io \cite{FG}, and Munteanu, Sesum \cite{MS} generalized the classification of locally conformally flat gradient shrinking Ricci solitons by showing that  a gradient shrinking Ricci solitons with harmonic Weyl curvature is either Einstein, or a finite quotient of $N^k\times\mathbb{R}^{n-k}$ for $0\leq k\leq n$, where $N^k$ is a $k$-dimensional Einstein manifold of positive scalar curvature. H.-D. Cao, Q. Chen \cite{CaoChen} also proved that a Bach-flat gradient shrinking Ricci soliton is either Einstein, or a finite quotient of $N^{n-1}\times\mathbb{R}$ or $\mathbb{R}^n$, where $N^{n-1}$ is an $(n-1)$-dimensional Einstein manifold. Naber \cite{Naber}, Munteanu, M.-T. Wang \cite{MW}, Catino \cite{Catino}, and M. Cai \cite{Cai} proved various rigidity results under appropriate curvature pinching assumptions.

Dimension four is the lowest dimension where there are interesting examples of shrinking gradient Ricci solitons. The first examples where constructed by Koiso \cite{Koiso} and Cao \cite{Cao2}, also see \cite{FIK, WZ}. Note that all of the known interesting examples are K\"ahler. In dimension 4, the Hodge star splits the space of $2$-forms into the self-dual and anti-self dual parts and consequently the curvature tensor and Weyl tensor respect this splitting (see section 2).  It is thus natural to consider self dual or anti-self dual part of Weyl curvature $W^{\pm}$ commonly called the half Weyl curvature. X. Chen, Y. Wang \cite{ChenWang} and H.-D. Cao, Q. Chen \cite{CaoChen} proved that a half conformally flat ($W^{\pm}=0$) four-dimensional gradient shrinking Ricci soliton is a finite quotient of $S^4$, $\mathbb{C}P^2$, $S^3\times\mathbb{R}$, or $\mathbb{R}^4$. In \cite{Wu2}, the second author proved that a compact four-dimensional gradient shrinking Ricci soliton with $\delta W^{\pm}=0$ and half two-nonnegative curvature operator (which is equivalent to half nonnegative isotropic curvature) is a finite quotient of $S^4$ or K\"ahler-Einstein.

In this paper we complete the classification of four-dimensional gradient shrinking Ricci solitons with harmonic half Weyl curvature,

\begin{thm} \label{Thm1.1}
A four-dimensional gradient shrinking Ricci soliton with $\delta W^{\pm}=0$ is either Einstein, or a finite quotient of $S^3\times\mathbb{R}$, $S^2\times\mathbb{R}^2$, or $\mathbb{R}^4$.
\end{thm}

As opposed to the locally conformally flat or half-locally conformally flat conditions we note that all Einstein metrics satisfy $\delta W = 0$  and are Bach flat and all Einstein $4$-manifolds satisfy $\delta W ^{\pm} = 0$. Thus, for example Theorem \ref{Thm1.1} shows that one can differentiate whether a gradient Ricci soliton is Einstein or non-Einstein simply by checking $\delta W^{\pm}$ (or the Bach tensor).  There are many interesting results about the tensor $W^{\pm}$ in the geometry of $4$-manifolds see for example Gursky \cite{Gursky} for an interesting gap theorem for $\|W^{\pm}\|_{L^2}$.

We also note that any $4$-dimensional K\"ahler metric with constant scalar curvature satisfies $\delta W^+=0$ and $\frac{|W^+|^2}{R^2}=\frac{1}{24}$ (see \cite{Der}). In the process of  detecting $S^2\times\mathbb{R}^2$ in the proof of Theorem \ref{Thm1.1}, we prove,
\begin{thm} \label{Thm1.3}
A four-dimensional gradient K\"ahler-Ricci soliton with constant scalar curvature is either K\"ahler-Einstein, or a finite quotient of $M\times\mathbb{C}$, where $M$ is a Riemann surface.
\end{thm}

Theorem \ref{Thm1.3} has also been obtained by Fernandez-Lopez and Garcia-Rio \cite{FG2} recently using a different method, and they were also able to classify six-dimensional gradient K\"ahler-Ricci solitons with constant scalar curvature.

Recently there have been some other interesting classification results for K\"ahler-Ricci solitons. Y. Su and K. Zhang \cite{SZ} have proved that a complete noncompact gradient K\"ahler-Ricci soliton with vanishing Bochner tensor is K\"ahler-Einstein. Q. Chen, M. Zhu \cite{CZ} have proved that a gradient steady K\"ahler-Ricci soliton with harmonic Bochner tensor is Calabi-Yau, and a gradient shrinking (expanding) K\"ahler-Ricci soliton with harmonic Bochner tensor is either K\"ahler-Einstein, or a finite quotient of $N^k\times\mathbb{C}^{n-k}$, where $N^k$ is a K\"ahler-Einstein manifold of positive (negative) scalar curvature. Calamai, Petrecca \cite{CP} also have proved that an extremal K\"ahler-Ricci soliton with positive holomorphic sectional curvature is K\"ahler-Einstein.

In the process of proving  Theorem \ref{Thm1.1} we also observe that in a result of Catino \cite{Catino}, the nonnegative Ricci curvature assumption is not required (by replacing $|\mathrm{Ric}|^2$
by $|\mathrm{Ric}|$ in his argumnet),
\begin{thm}[Catino \cite{Catino}] \label{Thm1.5}
A gradient shrinking Ricci soliton with
\begin{equation*}
|W|R\leq\sqrt{\frac{2(n-1)}{n-2}}\left(\overset{\circ}{\mathrm{Ric}}-\frac{R}{\sqrt{n(n-1)}}\right)^2
\end{equation*}
is a finite quotient of $S^n$, $S^{n-1}\times\mathbb{R}$, or $\mathbb{R}^n$.
\end{thm}

\begin{remark}
We also point out that our proof is different from the $\delta W=0$ case. For gradient shrinking Ricci solitons with $\delta W=0$, the proofs of Fern\'andez-L\'opez, Garc\'ia-R\'io \cite{FG}, and Munateanu, Sesum \cite{MS} rely on the following identity. If the Ricci curvature is bounded below and $|\mathrm{Rm}|\leq e^{a(r+1)}$ for some $a\in\mathbb{R}$ (X. Cao \cite{CaoX}), or if $\int_M|\mathrm{Rm}|^2e^{-\delta f}<\infty$ for some $\delta<1$ (Munteanu, Sesum \cite{MS}), then
\begin{equation*}
\int_M |\delta\mathrm{Rm}|^2e^{-f}dv=\int_M|\nabla\mathrm{Ric}|^2e^{-f}dv.
\end{equation*}
Unfortunately, it is not clear whether there is an analogous identity for half curvature.
\end{remark}

The main arguments to prove Theorem \ref{Thm1.1} are curvature decompostions, the Weitzenb\"ock formula for half Weyl curvature \cite{CT,Wu2}, and an analogous argument of Catino \cite{Catino}. As with much of the work mentioned above, another important component is the $D$-tensor introduced by Cao and Chen in \cite{CaoChen2, CaoChen}. We observe that the $D$-tensor arises naturally from the standard curvature decomposition. First we show, using the decomposition  of the curvature tensor, that if $\delta W^{\pm}=0$ then $\nabla f$ is an eigenvector of the Ricci tensor, and observe that $W^{\pm}$ can be expressed explicitly in terms of the traceless Ricci curvature. Next applying the Weitzenb\"ock formula and the maximum principle, we prove 
an identity involving the (anti-)self-dual Weyl curvature, the traceless Ricci curvature and the scalar curvature, which further imply that either $\nabla f\equiv 0$, hence $(M,g)$ is an Einstein manifold; or $W^{\pm}\equiv 0$, hence $(M,g)$ is a finite quotient of $S^3\times\mathbb{R}$ or $\mathbb{R}^4$; or $(M,g)$ is a cscK gradient shrinking K\"ahler-Ricci soliton, hence is a finite quotient of $S^2\times\mathbb{R}^2$.

\

The rest of the paper is organized as follows. In section 2, we discuss curvature decompositions and the relationship between $W^{\pm}$ and the traceless Ricci curvature when $\delta W^{\pm}=0$. In section 3, we prove Theorem \ref{Thm1.1}, and Theorem \ref{Thm1.3} by applying the Weitzenb\"ock formula and the maximum principle.

\textbf{Acknowledgement}.
Part of the work was done when the first author was visiting the Department of Mathematics at Cornell University, he greatly thanks Professor Xiaodong Cao for his help and Department of Mathematics for their hospitality. The second author thanks Professors Jeffrey Case, Yuanqi Wang and Yuan Yuan for helpful discussions. The first author is partially supported by NSFC (11101267, 11271132) and the China Scholarship Council (201208310431). The second author is partially supported by an AMS-Simons postdoc travel grant.


\section{Curvature decompositions on four-dimensional gradient Ricci solitons}

In this section we will discuss curvature decompositions on four-dimensional gradient Ricci solitons. First we fix some notation.  Our sign conventions for the curvature tensor will be so that

\begin{equation*}
\begin{split}
&R_{ijkl}=g_{hk}R^h_{ijl},\ \ K(e_i,e_j)=R_{ijij},\ \
R_{ik}=g^{jl}R_{ijkl},\ \ R=g^{ij}R_{ij}.\\
\end{split}
\end{equation*}
And our convention for the inner product of two $(0,4)$-tensors $S,T$ will be
\begin{equation*}
\langle S,T\rangle=\frac{1}{4}S_{ijkl}T^{ijkl}\
\end{equation*}
so that our convention agrees with the  one in Derdzinski's Weitzenb\"ock formula \cite{Der}.

On an oriented $4$-manifold, the Hodge star $\star: \wedge^2 M \rightarrow \wedge^2 M$ has eigenvalues $1$ and $-1$.  Thus we can break $\wedge^2 M = \wedge^+ M + \wedge^-M$ according to the eigenspaces of $\star$.  Given a basis $\{ e_1, e_2, e_3, e_4 \}$ of $T_pM$,  for any pair $(ij)$, $1\leq i\neq j\leq 4$, denote $(i'j')$ to be the dual of $(ij)$, i.e., the pair such that $e_i\wedge e_j\pm e_{i'}\wedge
e_{j'}\in \wedge^{\pm}M$. In other words, $(iji'j')=\sigma(1234)$ for some even permutation $\sigma\in S_4$.
So for any $(0,4)$-tensor $T$, its (anti-)self-dual part is
\begin{equation*}
T^{\pm}_{ijkl}=\frac{1}{4}(T_{ijkl}\pm T_{ijk'l'}\pm T_{i'j'kl}+T_{i'j'k'l'}),
\end{equation*}

It is well known that for four-manifolds, Weyl curvature has a very interesting symmetry,

\begin{lemma} \label{Lemma2.2}
Let $(M, g)$ be a four-dimensional Riemannian manifold. Then
\begin{equation*}
\begin{split}
W_{ijkl}=&W_{i'j'k'l'},
\end{split}
\end{equation*}
therefore,
\begin{equation*}
\begin{split}
W^{\pm}_{ijkl}=&\pm W^{\pm}_{ijk'l'}=\pm W^{\pm}_{i'j'kl}=W^{\pm}_{i'j'k'l'}=\frac{1}{2}(W_{ijkl}{\pm}W_{ijk'l'}).
\end{split}
\end{equation*}
In particular, for any $u\in C^{\infty}(M)$,
\begin{equation*}
\begin{split}
|\iota_{\nabla u}W^{\pm}|^2=&\frac{1}{4}|W^{\pm}|^2|\nabla u|^2.
\end{split}
\end{equation*}
\end{lemma}

\

Now we discuss curvature decompositions on four-dimensional gradient Ricci solitons. First recall that, similar to Einstein metrics, for gradient Ricci solitons, we have
\begin{equation*}
\begin{split}
\delta(e^{-f}\mathrm{Rm})=0,
\end{split}
\end{equation*}
following from which we have the following basic identities,

\begin{lemma} \label{Lemma2.1}
Let $(M, g, f)$ be a gradient Ricci soliton. Then
\begin{equation*}
\begin{split}
\nabla_kR_{jl}-\nabla_lR_{jk}=&R_{ijkl}\nabla^if,\\
(\delta\mathrm{Rm})_{jkl}=\nabla^iR_{ijkl}=&R_{ijkl}\nabla^if,\\
\nabla_iR=2\nabla^jR_{ij}=&2R_{ij}\nabla^jf.
\end{split}
\end{equation*}
\end{lemma}

In \cite{CaoChen2, CaoChen}, H.-D. Cao, Q. Chen introduced a $(0,3)$-tensor $D=\frac{n-2}{n-3}\delta W-\iota_{\nabla f}W$, which plays an important role in their classification of locally conformally flat gradient steady Ricci solitons and Bach-flat gradient shrinking Ricci solitons. We observe that $D$ and its ``self-dual" and ``anti-self-dual" parts $D^{\pm}$, arise naturally from the standard curvature decomposition. For our purpose we only calculate the four-dimensional case, for general dimensions the argument is the same.

\begin{lemma} \label{Lemma2.4}
Let $(M,g,f)$ be a four-dimensional gradient Ricci soliton. Then
\begin{equation*}
\begin{split}
D_{jkl}=&2\nabla^iW_{ijkl}-W_{ijkl}\nabla^if\\
=&\frac{1}{2}(R_{jl}\nabla_kf-R_{jk}\nabla_lf)+\frac{1}{12}(\nabla_kRg_{jl}-\nabla_lRg_{jk})
-\frac{R}{6}(g_{jl}\nabla_kf-g_{jk}\nabla_lf),\\
D^{\pm}_{jkl}\overset{\triangle}{=}&2\nabla^iW^{\pm}_{ijkl}-W^{\pm}_{ijkl}\nabla^if\\
=&\frac{1}{4}(R_{jl}\nabla_kf-R_{jk}\nabla_lf)+\frac{1}{24}(\nabla_kRg_{jl}-\nabla_lRg_{jk})\\
-&\frac{R}{12}(g_{jl}\nabla_kf-g_{jk}\nabla_lf)\\
\pm&\frac{1}{4}(R_{jl'}\nabla_{k'}f-R_{jk'}\nabla_{l'}f)\pm\frac{1}{24}(\nabla_{k'}Rg_{jl'}-\nabla_{l'}Rg_{jk'})\\
\mp&\frac{R}{12}(g_{jl'}\nabla_{k'}f-g_{jk'}\nabla_{l'}f).
\end{split}
\end{equation*}
In particular,
\begin{equation}\label{Dnorm}
\begin{split}
|D^+|^2=|D^-|^2=\frac{1}{2}|D|^2=&\frac{1}{4}|\overset{\circ}{\mathrm{Rc}}|^2|\nabla f|^2-\frac{1}{48}|R\nabla f-2\nabla R|^2.
\end{split}
\end{equation}
\end{lemma}

\Pf. We apply the standard curvature decomposition to both sides of the identity $\nabla^iR_{ijkl}=R_{ijkl}\nabla^if$,

For the left hand side, we have 
\begin{equation*} \label{curvdecom1}
\begin{split}
\nabla^iR_{ijkl}=&\nabla^iW_{ijkl}+\frac{1}{2}\nabla^i(R_{ik}g_{jl}+R_{jl}g_{ik}-R_{il}g_{jk}-R_{jk}g_{il})\\
-&\frac{\nabla^iR}{6}(g_{ik}g_{jl}-g_{il}g_{jk})\\
=&\nabla^iW_{ijkl}+\frac{1}{2}(\nabla_kR_{jl}-\nabla_lR_{jk}+\frac{1}{2}\nabla_kRg_{jl}-\frac{1}{2}\nabla_lRg_{jk})\\
-&\frac{1}{6}(\nabla_kRg_{jl}-\nabla_lRg_{jk})\\
=&\nabla^iW_{ijkl}+\frac{1}{2}\nabla^iR_{ijkl}+\frac{1}{4}(\nabla_kRg_{jl}-\nabla_lRg_{jk})
-\frac{1}{6}(\nabla_kRg_{jl}-\nabla_lRg_{jk}),
\end{split}
\end{equation*}
Therefore we get
\begin{equation*}
\begin{split}
\nabla^iR_{ijkl}=&2\nabla^iW_{ijkl}+\frac{1}{6}(\nabla_kRg_{jl}-\nabla_lRg_{jk}),
\end{split}
\end{equation*}

On the other hand, observe that by the second Bianchi identity, we have $\nabla^iR_{i'j'kl}=0$, so we get
\begin{equation*}
\begin{split}
4\nabla^iR^{\pm}_{ijkl}=&\nabla^i(R_{ijkl}+R_{ijk'l'}+R_{i'j'kl}+R_{i'j'k'l'})\\
=&\nabla^i(R_{ijkl}+R_{ijk'l'})\\
=&2\nabla^i(W_{ijkl}+W_{ijk'l'})+\frac{1}{6}(\nabla_kRg_{jl}-\nabla_lRg_{jk})
+\frac{1}{6}(\nabla_{k'}Rg_{jl'}-\nabla_{l'}Rg_{jk'})\\
=&4\nabla^iW^{\pm}_{ijkl}+\frac{1}{6}(\nabla_kRg_{jl}-\nabla_lRg_{jk})+\frac{1}{6}(\nabla_{k'}Rg_{jl'}-\nabla_{l'}Rg_{jk'}),
\end{split}
\end{equation*}
therefore we obtain
\begin{equation} \label{deltaRm}
\begin{split}
R_{ijkl}\nabla^if+R_{ijk'l'}\nabla^if=&4\nabla^iW^{\pm}_{ijkl}+\frac{1}{6}(\nabla_kRg_{jl}-\nabla_lRg_{jk})\\
+&\frac{1}{6}(\nabla_{k'}Rg_{jl'}-\nabla_{l'}Rg_{jk'})
\end{split}
\end{equation}

For the right hand side, we have
\begin{equation*} \label{curvdecom}
\begin{split}
R_{ijkl}\nabla^if=&W_{ijkl}\nabla^if+\frac{1}{2}(R_{ik}g_{jl}+R_{jl}g_{ik}-R_{il}g_{jk}-R_{jk}g_{il})\nabla^if\\
-&\frac{R}{6}(g_{ik}g_{jl}-g_{il}g_{jk})\nabla^if\\
=&W_{ijkl}\nabla^if+\frac{1}{2}(R_{jl}\nabla_kf-R_{jk}\nabla_lf+\frac{1}{2}\nabla_kRg_{jl}-\frac{1}{2}\nabla_lRg_{jk})\\
-&\frac{R}{6}(g_{jl}\nabla_kf-g_{jk}\nabla_lf).
\end{split}
\end{equation*}
Taking the difference, we get
\begin{equation*}
\begin{split}
D_{jkl}=&2\nabla^iW_{ijkl}-W_{ijkl}\nabla^if\\
=&\frac{1}{2}(R_{jl}\nabla_kf-R_{jk}\nabla_lf)+\frac{1}{12}(\nabla_kRg_{jl}-\nabla_lRg_{jk})
-\frac{R}{6}(g_{jl}\nabla_kf-g_{jk}\nabla_lf),
\end{split}
\end{equation*}
Therefore
\begin{equation*}
\begin{split}
D^{\pm}_{jkl}=&2\nabla^iW^{\pm}_{ijkl}-W^{\pm}_{ijkl}\nabla^if\\
=&\frac{1}{2}[(2\nabla^iW_{ijkl}-W_{ijkl}\nabla^if)\pm(2\nabla^iW_{ijk'l'}-W_{ijk'l'}\nabla^if)]\\
=&\frac{1}{4}(R_{jl}\nabla_kf-R_{jk}\nabla_lf)+\frac{1}{24}(\nabla_kRg_{jl}-\nabla_lRg_{jk})
-\frac{R}{12}(g_{jl}\nabla_kf-g_{jk}\nabla_lf)\\
\pm&\frac{1}{4}(R_{jl'}\nabla_{k'}f-R_{jk'}\nabla_{l'}f)\pm\frac{1}{24}(\nabla_{k'}Rg_{jl'}-\nabla_{l'}Rg_{jk'})\\
\mp&\frac{R}{12}(g_{jl'}\nabla_{k'}f-g_{jk'}\nabla_{l'}f),
\end{split}
\end{equation*}
\qed

Fern\'andez-L\'opez and Garc\'ia-R\'io \cite{FG} proved that if a gradient Ricci soliton satisfies $\delta W=0$, then $\nabla f$ is an eigenvector of the Ricci tensor. Following from Lemme \ref{Lemma2.4}, it is easy to see that in dimensions four, $\delta W^{\pm}=0$ provides the same information, see also X. Cao, Tran \cite{CT},

\begin{lemma} \label{Lemma2.3}
Let $(M,g,f)$ be a four-dimensional gradient Ricci soliton. If $\delta W^{\pm}=0$, then $\nabla f$, whenever nonzero, is an eigenvector of the Ricci tensor.
\end{lemma}

\Pf. In equation \eqref{deltaRm}, if $\delta W^{\pm}=0$, then
\begin{equation*}
\begin{split}
R_{ijkl}\nabla^if+R_{ijk'l'}\nabla^if=&\frac{1}{6}(\nabla_kRg_{jl}-\nabla_lRg_{jk})+\frac{1}{6}(\nabla_{k'}Rg_{jl'}-\nabla_{l'}Rg_{jk'})
\end{split}
\end{equation*}

Let $e_1=\frac{\nabla f}{|\nabla f|}$, and extend it to an orthonormal basis $\{e_1, e_2, e_3, e_4\}$ of $T_pM$.
Let $j=k=1$, $l\neq 1$, then since $(klk'l')=\sigma(1234)$, we have $g_{jl}=g_{jk'}=g_{jl'}=0$. Therefore we get
\begin{equation*}
\begin{split}
0=&\nabla_lR=2R_{lj}\nabla^jf=2|\nabla f|R_{1l}.
\end{split}
\end{equation*}
\qed

Combining Lemma \ref{Lemma2.4} and Lemma \ref{Lemma2.3}, we make a key observation that if $\delta W^{\pm}=0$, then $W^{\pm}$ has a nice expression in terms of
Ricci curvature and scalar curvature, 

\begin{prop} \label{Prop2.3}
Let $(M,g,f)$ be a four-dimensional gradient Ricci soliton with $\delta W^{\pm}=0$. Denote $a_1,a_2,a_3,a_4$ be the eigenvalues of the traceless Ricci tensor 
with corresponding eigenvectors $e_1=\frac{\nabla f}{|\nabla f|}, e_2, e_3, e_4$. Then whenever $\nabla f\neq 0$,
\begin{equation*}
\begin{split}
b_1\overset{\triangle}{=}W^{\pm}_{1212}=&-\frac{1}{12}(a_1+3a_2)=\frac{1}{12}(a_3+a_4-2a_2),\\
b_2\overset{\triangle}{=}W^{\pm}_{1313}=&-\frac{1}{12}(a_1+3a_3)=\frac{1}{12}(a_2+a_4-2a_3),\\
b_3\overset{\triangle}{=}W^{\pm}_{1414}=&-\frac{1}{12}(a_1+3a_4)=\frac{1}{12}(a_2+a_3-2a_4),\\
W^{\pm}_{1j1l}=&0,\quad \mathrm{ if }\ j\neq l.
\end{split}
\end{equation*}
\end{prop}

\Pf. By Lemma \ref{Lemma2.3} we have   $R_{1j}=0$ for $j\neq 1$, which gives us
\begin{equation*}
\nabla_1R=2R_{1j}\nabla^jf=2R_{11}|\nabla f|.
\end{equation*}
If $\delta W^{\pm}=0$, then
\begin{equation*}
\begin{split}
-W^{\pm}_{1j1j}|\nabla f|=&-W^{\pm}_{ij1j}\nabla^if\\
=&\frac{1}{4}(R_{jj}\nabla_1f-R_{1j}\nabla_jf)+\frac{1}{24}(\nabla_1Rg_{jj}-\nabla_jRg_{1j})\\
-&\frac{R}{12}(g_{jj}\nabla_1f-g_{1j}\nabla_jf)\\
\pm&\frac{1}{4}(R_{jj'}\nabla_{1'}f-R_{j1'}\nabla_{j'}f)\pm\frac{1}{24}(\nabla_{1'}Rg_{jj'}-\nabla_{j'}Rg_{j1'})\\
\mp&\frac{R}{12}(g_{jj'}\nabla_{1'}f-g_{j1'}\nabla_{j'}f)\\
=&\frac{1}{4}R_{jj}|\nabla f|+\frac{1}{12}R_{11}|\nabla f|-\frac{R}{12}|\nabla f|.
\end{split}
\end{equation*}
If $\nabla f\neq 0$, then we get
\begin{equation*}
\begin{split}
-W^{\pm}_{1212}=&\frac{1}{4}R_{22}+\frac{1}{12}R_{11}-\frac{R}{12}\\
=&\frac{1}{12}[3(R_{22}-\frac{R}{4})+(R_{11}-\frac{R}{4})]\\
=&\frac{1}{12}(a_1+3a_2),
\end{split}
\end{equation*}
similarly we get $W^{\pm}_{1313}$ and $W^{\pm}_{1414}$.

If $j\neq l$, then it is easy to compute that $W^{\pm}_{1j1l}=0$.

\qed


\section{Proof of Theorem 1.1}

First recall the Weitzenb\"ock formula for $W^{\pm}$ (see \cite{CT} or \cite{Wu2}),

\begin{prop} \label{Prop3.1}
Let $(M, g, f)$ be a four-dimensional gradient Ricci soliton. Then
\begin{equation*}
\begin{split}
\Delta_f|W^{\pm}|^2=&2|\nabla W^{\pm}|^2+4\lambda|W^{\pm}|^2-36\det
W^{\pm}-\langle(\overset{\circ}{\mathrm{Ric}}\circ\overset{\circ}{\mathrm{Ric}})^{\pm},W^{\pm}\rangle.
\end{split}
\end{equation*}
\end{prop}

\begin{remark}
The second author \cite{Wu2} derived a Weitzenb\"ock formula for generalized $m$-quasi-Einstein metrics (or ``Einstein metrics" on smooth metric measure spaces), which are defined by
\begin{equation*}
\begin{split}
\mathrm{Ric}+\nabla^2f-\frac{1}{m}df\otimes df=\lambda g,
\end{split}
\end{equation*}
for some $f,\lambda\in C^{\infty}(M)$ and an $m\in\mathbb{R}\cup\{\pm\infty\}$, from which we expect a similar rigidity result for quasi-Einstein four-manifolds with $\delta W^{\pm}=0$.
\end{remark}

Next we compute,

\begin{prop} \label{Prop3.2}
Let $(M, g, f)$ be a four-dimensional gradient Ricci soliton. Let $h=f-\ln R^2$, then
\begin{equation} \label{quotient}
\begin{split}
\Delta_h\left(\frac{|W^{\pm}|}{R}\right)\geq&
\frac{1}{2|W^{\pm}|R^2}\Big(R^2|W^{\pm}|^2-36R\det W^{\pm}
+4|W^{\pm}|^2|\overset{\circ}{\mathrm{Ric}}|^2\\
&-R\langle(\overset{\circ}{\mathrm{Ric}}\circ\overset{\circ}{\mathrm{Ric}})^{\pm},W^{\pm}\rangle\Big).
\end{split}
\end{equation}
\end{prop}

\Pf. Recall the Kato inequality $|\nabla T|^2\geq|\nabla|T||^2$ for any tensor $T$. From Proposition \ref{Prop3.1}, we get
\begin{equation*}
\begin{split}
\Delta_f|W^{\pm}|=&\frac{1}{2|W^{\pm}|}\Big[2|\nabla W^{\pm}|^2-2|\nabla|W^{\pm}||^2+4\lambda|W^{\pm}|^2-36\det W^{\pm}\\
&-\langle(\overset{\circ}{\mathrm{Ric}}\circ\overset{\circ}{\mathrm{Ric}})^{\pm},W^{\pm}\rangle\Big]\\
\geq&\frac{1}{2|W^{\pm}|}\Big[4\lambda|W^{\pm}|^2-36\det W^{\pm} -\langle(\overset{\circ}{\mathrm{Ric}}\circ\overset{\circ}{\mathrm{Ric}})^{\pm},W^{\pm}\rangle\Big]
\end{split}
\end{equation*}
and recall that
\begin{equation} \label{Requation}
\begin{split}
\Delta_fR=&2\lambda R-2|\mathrm{Ric}|^2.
\end{split}
\end{equation}
Therefore we compute
\begin{equation*}
\begin{split}
\Delta_f\left(\frac{|W^{\pm}|}{R}\right)
=&\frac{\Delta_f|W^{\pm}|}{R}-\frac{|W^{\pm}|\Delta_f R}{R^2}-2\frac{\nabla|W^{\pm}|\nabla R}{R^2}+2\frac{|W^{\pm}||\nabla R|^2}{R^3}\\
\geq&\frac{1}{2|W^{\pm}|}(4\lambda|W^{\pm}|^2-36\det W^{\pm}
-\langle(\overset{\circ}{\mathrm{Ric}}\circ\overset{\circ}{\mathrm{Ric}})^{\pm},W^{\pm}\rangle)\\
-&\frac{|W^{\pm}|}{R^2}(2\lambda R-2|\mathrm{Ric}|^2)\\
-&2\frac{\nabla|W^{\pm}|\nabla R}{R^2}+2\frac{|W^{\pm}||\nabla R|^2}{R^3}\\
=&-2\frac{1}{R}\Big\langle\nabla\Big(\frac{|W^{\pm}|}{R}\Big),\nabla R\Big\rangle
+\frac{1}{2|W^{\pm}|R^2}\Big(4|W^{\pm}|^2|\mathrm{Ric}|^2\\
-&36R\det W^{\pm}-R\langle(\overset{\circ}{\mathrm{Ric}}\circ\overset{\circ}{\mathrm{Ric}})^{\pm},W^{\pm}\rangle\Big)\\
=&-\Big\langle\nabla\Big(\frac{|W^{\pm}|}{R}\Big),\nabla\ln R^2\Big\rangle+\frac{1}{2|W^{\pm}|R^2}\Big(R^2|W^{\pm}|^2\\
-&36R\det W^{\pm}+4|W^{\pm}|^2|\overset{\circ}{\mathrm{Ric}}|^2
-R\langle(\overset{\circ}{\mathrm{Ric}}\circ\overset{\circ}{\mathrm{Ric}})^{\pm},W^{\pm}\rangle\Big).
\end{split}
\end{equation*}

\qed

We have,

\begin{lemma} \label{Lemma3.1}
Let $(M, g, f)$ be a four-dimensional gradient shrinking Ricci soliton with $\delta W^{\pm}=0$, then whenever $\nabla f\neq 0$,
\begin{equation*}
\begin{split}
R^2|W^{\pm}|^2-36R\det W^{\pm}+4|W^{\pm}|^2|\overset{\circ}{\mathrm{Ric}}|^2
-R\langle(\overset{\circ}{\mathrm{Ric}}\circ\overset{\circ}{\mathrm{Ric}})^{\pm},W^{\pm}\rangle\geq 0,
\end{split}
\end{equation*}
with equality if and only if, either

$\mathrm{(1)}$, $a_2=a_3=a_4$, i.e., $W^{\pm}=0$; or

$\mathrm{(2)}$, $a_1=a_i=-a, a_j=a_k=a$, $2\leq i,j,k\leq 4$, and $R=4a$ for some $a>0$, where $a_1,a_2,a_3,a_4$ are eigenvalues of $\overset{\circ}{\mathrm{Ric}}$ with corresponding orthonormal eigenvectors $e_1=\frac{\nabla f}{|\nabla f|}, e_2, e_3, e_4$.
\end{lemma}

\begin{remark}
Lemma \ref{Lemma3.1} also works for gradient steady and expanding solitons, and the sign of $a$ in the second equality case changes correspondingly.
\end{remark}

The proof of the Lemma will be presented at the end of this section. We first prove Theorem \ref{Thm1.1}.

\

\Pf\ of Theorem \ref{Thm1.1}. B.-L. Chen \cite{Chen} proved that any gradient shrinking
Ricci soliton has $R\geq 0$. Moreover, either $R>0$ on $M$, or $R\equiv 0$ on $M$, and if $R\equiv 0$ then $(M,g)$
is a finite quotient of $\mathbb{R}^4$, see \cite{PW2, PRS}. From now on we assume $R>0$.

If $M$ is compact, then from Proposition \ref{Prop3.2},
\begin{equation*}
\begin{split}
0=&\int_M\Delta_h\left(\frac{|W^{\pm}|}{R}\right)e^{-h}dv\\
\geq&\int_M\frac{1}{2|W^{\pm}|}\Big(R^2|W^{\pm}|^2-36R\det W^{\pm}
+4|W^{\pm}|^2|\overset{\circ}{\mathrm{Ric}}|^2\\
&-R\langle(\overset{\circ}{\mathrm{Ric}}\circ\overset{\circ}{\mathrm{Ric}})^{\pm},W^{\pm}\rangle\Big)e^{-f}dv
\end{split}
\end{equation*}
so by Lemma \ref{Lemma3.1} we get
\begin{equation} \label{zero}
\begin{split}
&R^2|W^{\pm}|^2-36R\det W^{\pm}+4|W^{\pm}|^2|\overset{\circ}{\mathrm{Ric}}|^2
-R\langle(\overset{\circ}{\mathrm{Ric}}\circ\overset{\circ}{\mathrm{Ric}})^{\pm},W^{\pm}\rangle\equiv 0.
\end{split}
\end{equation}

If $M$ is noncompact, by equation \eqref{Dnorm} in Lemma \ref{Lemma2.4}, if $\delta W^{\pm}=0$, then
\begin{equation*}
\begin{split}
|W^{\pm}|\leq|\overset{\circ}{\mathrm{Ric}}|<|\mathrm{Ric}|,
\end{split}
\end{equation*}
Munteanu, Sesum \cite{MS} proved that for a gradient shrinking Ricci soliton,
\begin{equation*}
\begin{split}
\int_M |\mathrm{Ric}|^2e^{-\delta f}dv<\infty,
\end{split}
\end{equation*}
for any $\delta>0$. Therefore if $\delta W^{\pm}=0$, then
\begin{equation*}
\Big\|\frac{|W^{\pm}|}{R}\Big\|_{L^2_h(M)}=\int_M|W^{\pm}|^2e^{-f}dv<\infty.
\end{equation*}

By a maximum principle of Naber \cite{Naber} and Petersen and Wylie \cite{PW}, if $\int_M e^{-h}dv<\infty$, then any $L_h^2$-integrable $
h$-subharmonic function is constant, therefore we conclude that $\frac{|W^{\pm}|}{R}$ is constant, which also implies equation \eqref{zero}.

\

Recall that any gradient Ricci soliton is
a real-analytic manifold (see \cite{Ivey2} or \cite{Kotsch}), hence all $|\nabla f|^2$, $|W^{\pm}|^2$, and $R$
are analytic functions on $M$, therefore either $\nabla f\equiv 0$ or $W^{\pm}\equiv 0$, or the second equality
case in Lemma \ref{Lemma3.1} holds on $M$.

Case 1. If $\nabla f\equiv 0$ on $M$,then $(M,g)$ is Einstein.

Case 2. If $W^{\pm}\equiv 0$ on $M$, then $(M,g)$ is a finite quotient of $S^3\times\mathbb{R}$ or $\mathbb{R}^4$ by \cite{ChenWang, CaoChen}.

Case 3. If $\nabla f\not\equiv 0$ and $W^{\pm}\not\equiv 0$, then the second case in Lemma \ref{Lemma3.1} holds in an open dense set $S$ of $M$. Without loss of generality, assume $\delta W^+=0$ (if $\delta W^-=0$ changing the orientation we get $\delta W^+=0$).
First it is easy to see that
\begin{equation*}
R_{11}=\overset{\circ}{R}_{11}+\frac{R}{4}=0,\quad R_{ii}=0,\quad R_{jj}=R_{kk}=2a,
\end{equation*}
$2\leq i,j,k\leq 4$, therefore by Lemma \ref{Lemma2.3}, $\nabla R=\nabla_1Re_1=2R_{11}|\nabla f|=0$, that is $R\equiv$const on $S$. By the continuity of $R$, we have $R\equiv$const on $M$, and $|W^{\pm}|^2=\frac{R^2}{24}\equiv$const.
Moreover in this case, $W^{\pm}$ has only two distinct eigenvalues:
$$\pm\frac{R}{6}, \mp\frac{R}{12}, \mp\frac{R}{12},$$
(in fact $\mathfrak{R}^{\pm}$ is two-nonnegative), hence by a Theorem of Derdzinski (Proposition 5 in \cite{Der}), $g$ is a K\"ahler metric.

Assume $M$ is compact, since $R\equiv$const, by the soliton equation $R+\Delta f=4\lambda$, so $f\equiv$const (and $g$ is K\"ahler-Einstein), contradiction!
Therefore $(M,g,f)$ is a complete noncompact cscK gradient shrinking K\"ahler-Ricci soliton.

Furthermore, since the eigenvalues of Ricci curvature are $0,0,2a,2a$, so we have $|\mathrm{Ric}|^2=8a^2$, plugging into equation \eqref{Requation}, we get
\begin{equation*}
\begin{split}
0=\Delta_fR=&2\lambda R-2|\mathrm{Ric}|^2\\
=&8\lambda a-16a^2,
\end{split}
\end{equation*}
therefore $a=\frac{\lambda}{2}$, which in particular implies that $0\leq\mathrm{Ric}\leq \lambda g$. By Proposition 1.3 in \cite{PW2}, $(M,g,f)$ is rigid, i.e., it is a finite quotient of $N^k\times\mathbb{R}^{4-k}$, where $N^k$ is an Einstein manifold. Since $g$ is K\"ahler, so $(M,g)$ is a finite quotient of $N^2\times\mathbb{R}^2$. Since $R$ is positive, $N^2$ has to be $S^2$, therefore $(M,g)$ is a finite quotient of $S^2\times\mathbb{R}^2$.

\qed

\Pf\ of Theorem \ref{Thm1.3}. On a four-manifold, a K\"ahler metric with constant scalar curvature satisfies
\begin{equation*}
\delta W^+=0,\quad \frac{|W^+|^2}{R^2}=\frac{1}{24}.
\end{equation*}

If $(M,g,f)$ is a gradient shrinking Ricci soliton, it follows directly from Theorem \ref{Thm1.1}.

If $(M,g,f)$ is a gradient steady Ricci soliton, then $R\equiv$const implies that $R\equiv 0$ (see \cite{FG3} or \cite{Wu}), and by Proposition 4.3 in \cite{PW2}, $(M,g,f)$ is a finite quotient of $M\times\mathbb{C}$, where $M$ is a flat Riemann surface.

If $(M,g,f)$ is a gradient expanding Ricci soliton, then by Proposition \ref{Prop3.2} and Lemma \ref{Lemma3.1} we have
\begin{equation*}
\begin{split}
R^2|W^+|^2-36R\det W^+ +4|W^+|^2|\overset{\circ}{\mathrm{Ric}}|^2
-R\langle(\overset{\circ}{\mathrm{Ric}}\circ\overset{\circ}{\mathrm{Ric}})^+,W^+\rangle\equiv0.
\end{split}
\end{equation*}
Similar to the proof of Theorem \ref{Thm1.1}, there are three cases,

Case 1. If $\nabla f\equiv 0$ on $M$, then $(M,g)$ is K\"ahler-Einstein.

Case 2. If $W^+\equiv 0$ on $M$, then $D^+\equiv 0$, hence by Lemma \ref{Lemma2.4}, $D\equiv 0$, therefore $W\equiv0$ by Theorem 5.1 in H.-D. Cao and Q. Chen \cite{CaoChen}, and $(M,g,f)$ is a finite quotient of Gaussian expanding soliton by Y. Su and K. Zhang \cite{SZ}.

Case 3. If $\nabla f\not\equiv 0$ and $W^+\not\equiv 0$, then it follows from Case 3 in the proof of Theorem \ref{Thm1.1} that it is rigid, hence a finite quotient of $M\times\mathbb{C}$, where $M$ is a Riemann surface of constant negative curvature.

\qed

\Pf\ of Lemma \ref{Lemma3.1}. By proposition \ref{Prop2.3}, we express each term in terms of eigenvalues of traceless Ricci tensor,
\begin{equation*}
\begin{split}
|\overset{\circ}{\mathrm{Rc}}|^2=&a_1^2+a_2^2+a_3^2+a_4^2\\
=&2(a_2^2+a_3^2+a_4^2+a_2a_3+a_2a_4+a_3a_4),\\
|W^{\pm}|^2=&4(b_1^2+b_2^2+b_3^2)\\
=&\frac{1}{6}(a_2^2+a_3^2+a_4^2-a_2a_3-a_2a_4-a_3a_4)\\
36\det W^{\pm}=&\frac{1}{6}(a_3+a_4-2a_2)(a_2+a_4-2a_3)(a_2+a_3-2a_4)\\
=&\frac{1}{6}(-2a_2^3-2a_3^3-2a_4^3+3a_2^2a_3+3a_3^2a_2+3a_2^2a_4\\
&\hspace{1.5cm}+3a_4^2a_2+3a_3^2a_4+3a_4^2a_3-12a_2a_3a_4),\\
\langle(\overset{\circ}{\mathrm{Rc}}\circ\overset{\circ}{\mathrm{Rc}})^{\pm},W^{\pm}\rangle=&
W^{\pm}_{ijij}\overset{\circ}{R}_{ii}\overset{\circ}{R}_{jj}\\
=&2[b_1(a_1a_2+a_3a_4)+b_2(a_1a_3+a_2a_4)+b_3(a_1a_4+a_2a_3)]\\
=&\frac{1}{6}(2a_2^3+2a_3^3+2a_4^3+a_2^2a_3+a_3^2a_2+a_2^2a_4\\
&\hspace{1.5cm}+a_4^2a_2+a_3^2a_4+a_4^2a_3-12a_2a_3a_4).
\end{split}
\end{equation*}
Therefore we get
\begin{equation} \label{phi}
\begin{split}
6\phi=&6(R^2|W^{\pm}|^2-36R\det W^{\pm}+4|W^{\pm}|^2|\overset{\circ}{\mathrm{Ric}}|^2
-R\langle(\overset{\circ}{\mathrm{Ric}}\circ\overset{\circ}{\mathrm{Ric}})^{\pm},W^{\pm}\rangle)\\
=&R^2(a_2^2+a_3^2+a_4^2-a_2a_3-a_2a_4-a_3a_4)\\
-&4R(a_2^2a_3+a_3^2a_2+a_2^2a_4+a_4^2a_2+a_3^2a_4+a_4^2a_3-6a_2a_3a_4)\\
+&8(a_2^2+a_3^2+a_4^2+a_2a_3+a_2a_4+a_3a_4)\times\\
&\ \ (a_2^2+a_3^2+a_4^2-a_2a_3-a_2a_4-a_3a_4).
\end{split}
\end{equation}
By abusing the noataion, we identify $\phi$ and $6\phi$. Observe that $\phi$ is a fourth-order homogeneous symmetric polynomial
if we assume $R=k(a_2+a_3+a_4)$ for some $k\in\mathbb{R}$.

\

First we show $\phi\geq 0$ using Timofte's criterion for positivity of homogeneous symmetric polynomials (see Corollary 5.6 in \cite{Tim}),
\begin{prop}[Timofte \cite{Tim}]\label{Tim}
Let $p$ be a fourth-order homogeneous symmetric polynomial on $\mathbb{R}^n$, then
\begin{equation*}
p\geq 0\ \mathrm{on}\ \mathbb{R}^n \quad \Longleftrightarrow\quad
p(t\cdot\vec{1}_{\mathbb{R}^i},\vec{1}_{\mathbb{R}^{n-i}})\geq 0,\ \forall\ t\in[-1,1],\ i=1,2,...,n-1,
\end{equation*}
where $\vec{1}=(1,1,..,1)$.
\end{prop}

If $a_1\neq 0$, without of loss of generality, assume $R=-ka_1=k(a_2+a_3+a_4)$. In our case $n=3$, so we need to show that
\begin{equation*}
\phi(t,1,1)\geq 0,\quad \phi(t,t,1)\geq 0,\quad \forall\ t\in[-1,1].
\end{equation*}

For $\phi(t,1,1)$, plugging into equation \eqref{phi}, recall that $R=k(t+2)$, we get
\begin{equation*}
\begin{split}
\phi(t,1,1)=&(t-1)^2\Big[k^2(t+2)^2-8k(t+2)+8(t^2+2t+3)\Big].
\end{split}
\end{equation*}
Consider $\phi(t,1,1)$ as a quadratic function of $k$. When $-1\leq t\leq 1$, the discriminant
\begin{equation*}
\begin{split}
\mathfrak{D}=&-32(t+2)^2(t-1)^4(t+1)^2\leq 0,
\end{split}
\end{equation*}
and $\mathfrak{D}<0$ when $-1<t<1$. Therefore for all $-1\leq t\leq 1$,
\begin{equation*}
\begin{split}
\phi(t,1,1)\geq 0,
\end{split}
\end{equation*}
and $\phi(t,1,1)=0$ if and only if $t=1$, or $t=-1$ and $R=4$.

For $\phi(t,t,1)$, since we assume $a_1\neq 0$, so $t\neq-\frac{1}{2}$. Plugging into equation \eqref{phi}, recall that $R=k(2t+1)$, we get
\begin{equation*}
\begin{split}
\phi(t,t,1)=&(t-1)^2\Big[k^2(2t+1)^2-8kt(2t+1)+8(3t^2+2t+1)\Big].
\end{split}
\end{equation*}
Consider $\phi(t,t,1)$ as a quadratic function of $k$, we see that when $-1\leq t\leq 1$ and $t\neq-\frac{1}{2}$, the discriminant
\begin{equation*}
\begin{split}
\mathfrak{D}=&-32(2t+1)^2(t-1)^4(t+1)^2\leq 0,
\end{split}
\end{equation*}
and $\mathfrak{D}<0$ when $-1<t<1$ and $t\neq-\frac{1}{2}$. Therefore for all $-1\leq t\leq 1$ and $t\neq-\frac{1}{2}$,
\begin{equation*}
\begin{split}
\phi(t,t,1)\geq 0,
\end{split}
\end{equation*}
and $\phi(t,t,1)=0$ if and only if $t=1$, or $t=-1$ and $R=-4$.

\

If $a_1=0$ (which corresponds to $t=-\frac{1}{2}$ in $\phi(t,t,1)$), i.e., $a_2+a_3+a_4=0$, then $\phi$ can be simplified as
\begin{equation*}
\begin{split}
\phi=&3R^2(a_2^2+a_3^2+a_2a_3)-36Ra_2a_3(a_2+a_3)+24(a_2^2+a_3^2+a_2a_3)^2.
\end{split}
\end{equation*}
Consider $\phi$ as a quadratic function of $R$, then its discriminant
\begin{equation*}
\begin{split}
\mathfrak{D}=&36^2a_2^2a_3^2(a_2+a_3)^2-36[a_2^2+a_3^2+(a_2+a_3)^2]^3\leq 0.
\end{split}
\end{equation*}

Recall a well-known inequality: if $a+b+c=0$, then $3\sqrt{6}|abc|\leq(a^2+b^2+c^2)^{\frac{3}{2}}$,
with equality if and only if $a=-2b$, or $b=-2c$, or $c=-2a$. Since $-a_2-a_3+(a_2+a_3)=0$, we have
\begin{equation*}
\begin{split}
\mathfrak{D}=&36^2a_2^2a_3^2(a_2+a_3)^2-36[a_2^2+a_3^2+(a_2+a_3)^2]^3\\
\leq &24[a_2^2+a_3^2+(a_2+a_3)^2]^3-36[a_2^2+a_3^2+(a_2+a_3)^2]^3\\
\leq &0,
\end{split}
\end{equation*}
with equality if and only if $a_2=a_3=0$. So $\phi\geq 0$, and $\phi=0$ if and only if $a_2=a_3=a_4=0$.

Therefore we proved that $\phi(R,a_2,a_3,a_4)\geq 0$ on $\mathbb{R}^4$.

\

Next we show that when $a_2\neq a_3\neq a_4$, then $\phi>0$.

Assume that $a_2\neq a_3\neq a_4$ and $\phi(R,a_2,a_3,a_4)=0$. Taking the first derivatives we get,
\begin{equation*}
\begin{split}
\phi_{a_2}=&R^2(2a_2-a_3-a_4)-4R(2a_2a_3+2a_2a_4+a_3^2+a_4^2-6a_3a_4)\\
&\hspace{2.8cm}+16(2a_2^3+a_2a_3^2+a_2a_4^2-a_3^2a_4-a_3a_4^2-2a_2a_3a_4)=0,\\
\phi_{a_3}=&R^2(2a_3-a_2-a_4)-4R(2a_2a_3+2a_3a_4+a_2^2+a_4^2-6a_2a_4)\\
&\hspace{2.8cm}+16(2a_3^3+a_2^2a_3+a_3a_4^2-a_2^2a_4-a_2a_4^2-2a_2a_3a_4)=0,\\
\phi_{a_4}=&R^2(2a_4-a_2-a_2)-4R(2a_2a_4+2a_3a_4+a_2^2+a_3^2-6a_2a_3)\\
&\hspace{2.8cm}+16(2a_4^3+a_2^2a_4+a_3^2a_4-a_2^2a_3-a_2a_3^2-2a_2a_3a_4)=0.
\end{split}
\end{equation*}

Taking the difference, since $a_2\neq a_3\neq a_4$, we get
\begin{equation*}
\begin{split}
0=\frac{\phi_{a_2}-\phi_{a_3}}{a_2-a_3}=&3R^2-R(4(a_2+a_3)-32a_4)\\
&\hspace{0.6cm}+16(2a_2^2+2a_3^2+2a_4^2+a_2a_3+a_2a_4+a_3a_4),\\
0=\frac{\phi_{a_2}-\phi_{a_4}}{a_2-a_4}=&3R^2-R(4(a_2+a_4)-32a_3)\\
&\hspace{0.6cm}+16(2a_2^2+2a_3^2+2a_4^2+a_2a_3+a_2a_4+a_3a_4),\\
0=\frac{\phi_{a_3}-\phi_{a_4}}{a_3-a_4}=&3R^2-R(4(a_3+a_4)-32a_2)\\
&\hspace{0.6cm}+16(2a_2^2+2a_3^2+2a_4^2+a_2a_3+a_2a_4+a_3a_4).
\end{split}
\end{equation*}
Taking the difference again, we get $a_2=a_3=a_4$, contradiction! So $\phi>0$ when $a_2\neq a_3\neq a_4$.\\

Therefore $\phi\geq 0$, and by Timofte's criterion and above argument, $\phi=0$ if and only if,

either $a_2=a_3=a_4$, i.e. $W^{\pm}=0$;

or $a_1=a_i=-a, a_j=a_k=a$, $2\leq i,j,k\leq 4$, and $R=4a$, for some $a>0$.

\qed


\end{document}